\documentclass[11 pt]{amsart}
\usepackage{amsmath}
\usepackage{amssymb}
\usepackage{graphicx}

\newtheorem{theorem}{Theorem}[section]

\newtheorem{lemma}[theorem]{Lemma}
\newtheorem{proposition}[theorem]{Proposition}
\theoremstyle{definition}

\newtheorem{remark}[theorem]{Remark}

\numberwithin{equation}{section}
\usepackage{color,hyperref,tikz}
\hypersetup{urlcolor=blue, citecolor=blue}

\let\cal=\mathcal

\def\R{{\bf R}}

\def\mA{{\bf A}}

\def\mI{{\bf I}}

\def\skp#1#2#3{\langle\, #1 \mid #2 \,\rangle_{#3}}

\def\lD{{\cal D}}

\def\lK{{\cal K}}
\def\lN{{\cal N}}

\def\nor#1#2{{\| #1 \|}_{#2}}

\def\zi#1#2{[#1,#2]}

\def\eps{\varepsilon}

\def\diag{{\rm diag}}

\def\dv{{\rm div}}

\def\plav#1{\textcolor{black}{#1}}

\title[Greedy search of optimal approximate solutions]{Greedy search of optimal approximate solutions}
\author[M. Lazar]{Martin Lazar}    
	\address[M. Lazar]{Department of Electrical Engineering and Computing, University of Dubrovnik, 
\'Cira Cari\' ca 4, 20 000 Dub\-rov\-nik, Croatia }
\email{mlazar@unidu.hr}

\author[E. Zuazua]{ Enrique Zuazua}
\address[E. Zuazua]{Chair for Dynamics, Control and Numerics – Alexander von Humboldt Professorship, Department of Data Science. Friedrich-Alexander Universität, Erlangen-Nürnberg. 91058 Erlangen – Germany, \newline
Chair of Computational Mathematics, Deusto Foundation Av. de las Universidades 24, 48007 Bilbao, Basque Country, Spain\newline
Department of Mathematics, Autonomous University of Madrid, 28049 Mad\-rid, Spain \\
}
\email{enrique.zuazua@fau.de}
\date{\today}

\thanks{Dedicated to  Professor Ronald DeVore on the occasion of his 80th birthday.}


\begin{document}

\begin{abstract}
In this paper we develop a procedure to deal with a family of parameter-dependent ill-posed problems, for which the exact solution in general does not exist. The original problems are relaxed by considering corresponding approximate ones, whose optimal solutions are well defined, where the optimality is determined by the minimal norm requirement.
The procedure is based upon greedy algorithms that preserve, at least asymptotically, Kolmogorov approximation rates. In order to provide a-priori estimates for the algorithm, a Tychonoff-type regularization is applied, which adds an additional parameter to the model.
The theory is developed in an abstract theoretical framework that allows its application to  different  kinds of problems. We present a specific example that considers a family of ill-posed elliptic problems.
 The required general assumptions in this case translate to rather natural uniform lower and upper  bounds on coefficients of the considered operators.

\end{abstract}

\maketitle

\section{Introduction}

Greedy theory, inspired on the notion of Kolmogorov complexity and nonlinear approximation theory, has been extensively developed to provide optimal approximation rates for parameter depending problems, in particular Partial Differential Equations (PDE) \cite{VPRP03,CD15}.

This is a relevant subject since, often in applications, the model mimicking the dynamics is not fully known and is subject to uncertainty, in particular,  on some of the relevant parameters entering in the system such as diffusivity, Lamé coefficients etc.

A control theoretical counterpart was developed in \cite{LZ16}. The question addressed in that article referees to the classical problem of controllability, that of driving a dynamical system to a desired final configuration by the action of a suitable control. They did it in the context of parameter-depending problems, building a greedy algorithm allowing to guarantee an optimal approximation, in the sense of the Kolmogorov thickness, of the set of parameter-depending controls.

But, as observed in  \cite{LZ16}, an important case was left open. How can the greedy strategy be adapted when the system under consideration is not controllable, like it occurs for instance in parabolic PDEs in which, 
due to the strong time irreversibility of the model, only very smooth targets are reachable?

In this paper we consider these kind of problems in an abstract frame, motivated by that example, and characterized by the application of the greedy algorithm to linear systems in which the operators governing the system are not onto.

Inspired by the theory of approximate controllability for time-irreversible PDEs, developed by J. L. Lions, and by R. Glowinski and J. L. Lions in the numerical setting (e.g. \cite{L92,CGL94}), and the classical technique of Tychonoff regularization for ill-posed (inverse) problems, we adopt a two-folded perspective that allows us to end up developing a greedy strategy ensuring that Kolmogorov complexity is reached, for quasi-solutions (those that assure the fulfillment of the system up to an $\eps$ error) \plav{in those cases} where the exact solution  (corresponding to $\eps$ equal to zero) does not exist in the given functional setting.

With that purpose we proceed in several \plav{steps:}

1. We introduce a concept of optimal $\eps$-solution, reminiscent of the theory of approximate controllability, that allows characterizing the solution of minimal norm assuring that the system is solved up to an $\eps$ tolerance or error.

2. We link this $\eps$-solution with a suitable Tychonoff regularization, in which the $\eps$-error or tolerance, can be linked to the $\eta$-Tychonoff regularization parameter through a suitable nonlinear implicit equation.

3. We then adapt and apply the existing greedy methods for well posed parameter dependent problems, but in the context in which the number of free parameters is increased by one, to incorporate the $\eta$-Tychonoff parameter, that adds to the physical parameters on which the model depends.

4. We then use the nonlinear link of the $\eps$-error and $\eta$-Tychonoff parameter to derive a greedy approximation result for the $\eps$-problem.

Our abstract results apply to a wide variety of problems, such as deconvolution in image processing, time inversion of highly irreversible systems, like heat equations, or the approximate controllability of PDE.

Although we develop the theory in the context of bounded linear operators, it can also be extended to unbounded ones. This, in turn, allows us to study ill-posed elliptic problems in which the (exact) solution does in general  not exist for an arbitrary right-hand side. In this way, we generalize the setting considered by R. DeVore and his collaborators, in which they consider a family of well-posed elliptic problems and develop the greedy algorithms to recover the corresponding solutions \cite{CD15, DV17}.

The paper is organised as follows. The next section provides  definition  and analysis of optimal approximative solutions of a single $\eps$-problem. 
Section 3 is devoted to construction  of a greedy algorithm for solving  a family of parameter dependent problems. Generalisation and application of the theory to unbounded operators is presented in Section 4, followed by a particular example related to the elliptic equation. The paper is closed by some concluding remarks and perspectives for further research.

\section{Preliminaries on optimal approximate solutions of linear  systems }

Let $L$ be a bounded operator from $L(X,Y)$ with a dense range, where $X$ and $Y$ are infinite-dimensional, real Hilbert spaces. 
We consider the problem of solving the system 
\begin{equation}
  \label{sys_eq}
 Lu=f,
\end{equation}
for a given vector $f\in Y$.  Of course, due to the fact that image of  $L$ is only dense in $Y$, the above problem has no solution for every $f$. Thus we relax the problem and consider the approximate problem of finding $u$ such that 
\begin{equation}
  \label{appr_sol}
\nor{Lu-f}Y \leq \eps,
\end{equation}
for some a-priori given $ \eps >0$. 

Having assumed that the operator $L$ has a dense range, the set of approximate solutions satisfying \eqref{appr_sol} for any fixed $f$ is neither empty nor a singleton. 
Thus it is reasonable to choose the one which is optimal in some sense. We set the optimization criteria in this paper as the one of minimal norm. 
\plav{This leads us to the following problem
\begin{equation}
\label{opt_pr}
\min_{u\in X} \left\{ \frac{1}{2}\nor{u}X^2: \ \ \nor{Lu-f}Y \leq \eps \right\},
\end{equation}
whose solution we denote by $\tilde u$.
}

\plav{By exploring Fenchel-Rockafellar duality techniques  (e.g. \cite[\S  3.6]{Pey15}) the   solution $\tilde u$ of the problem \eqref{opt_pr} can be obtained by solving the corresponding dual problem.  This allows one to replace the original constrained problem by a non-constrained one. 
More precisely, the following theorem holds. 
\begin{proposition}
		\cite[Propositions 2 \& 4]{LM21}
		\label{GenHUM}
		The solution of problem $\eqref{opt_pr}$ is given by
		\begin{equation*}
		\tilde u=L^*  \tilde v,
		\label{form_sol}
		\end{equation*}
		where $\tilde v \in Y$ minimizes the dual problem
		\begin{equation}
  \label{funct_J}
J(v)=\frac{1}{2} \nor{L^* v}X^2+ \eps \nor{v}Y - \skp{f}{v}Y, \quad v\in Y.
\end{equation}
	\end{proposition}
	}

Here $L^*\in L(Y, X)$ stands for the adjoint operator of $L$, \plav{while $\skp{\cdot}{\cdot}Y$ denotes the scalar product in $Y$. }

\plav{
Note that the last proposition does not guarantee uniqueness of the solution to the dual problem \eqref{funct_J}. The latter includes the functional $J$ which}  is a non-standard one since it involves a non-smooth term of homogeneity one. 
Being strictly convex and continuous, it attains its minimal value at the unique point if it is coercive. 
According to the density assumption on $L$ and Hahn-Banach theorem, the adjoint operator $L^* $ is injective. But  in general it is not coercive (this would correspond to the case of $L$ having a full rank). However, the $\eps$ term entering the functional $J$ ensures its coercivity as shown by the following result.

\begin{proposition}
\label{J-coerc}
Functional $J$ defined by \eqref{funct_J} is coercive, i.e. it satisfies
\begin{equation*}
  \label{coerc}
  \liminf_{\nor{v}Y\to \infty} {J (v) \over \nor{v}Y} \geq \eps. 
   \end{equation*}
\end{proposition}
{\bf Proof:} 
The proof essentially follows the lines of the \cite[Proposition 2.1]{FPZ95} which treats the special case of $L$ being the heat operator.

We argue by contradiction. We suppose there exist a sequence $(v_n)$ such that $\nor{v_{n}}Y \to \infty$ and  
\begin{equation}
  \label{contr_coerc}
  \liminf_n {J (v_{n}) \over \nor{v_{n}}Y} < \eps. 
   \end{equation}
  Denoting by $v_n^0= v_n / \nor{v_{n}}Y$ the corresponding normalized vectors, it follows that, up to a subsequence, $(v_n^0)$ converges weakly to some $  v^0 \in Y$.
  From \eqref{contr_coerc} it follows
  $$
  \eps > \frac{1}{2}  \nor{v_{n}}Y \nor{L^* v_n^0}X^2+ \eps  - \skp{f}{v_n^0}Y.
$$
As all the terms in the last relation, except $\nor{v_{n}}Y$, are bounded, it implies that $\nor{L^* v_n^0}X \to 0 $. Consequently, $\nor{L^* v^0}X      =0$, and the injectivity of $L^* $ implies $v^0=0$.

Thus we obtain 
$$
 \liminf_n {J (v_{n}) \over \nor{v_{n}}Y} \geq \eps - \lim_n \skp{f}{v_n^0}Y =\eps,
 $$
 which contradicts the initial assumption \eqref{contr_coerc}.
\hfill $\Box$

\plav{
A detailed characterization of the unique minimizer of  $J$ is provided by the next result. 
\begin{proposition}
\label{prop_tilde_v}
The minimizer $\tilde v$ equals zero if and only if $\nor{f}Y \leq \eps$.  Otherwise, it
satisfies the    Euler-Lagrange equation which has the form
\begin{equation}
  \label{Eul-Lagr}
  L L^* \tilde v+ \eps {\tilde v \over \nor{\tilde v}Y} - f=0.
  \end{equation}
\end{proposition}
}
\plav{
{\bf Proof:} 
From the very definition  \eqref{funct_J} of functional $J$ it follows 
$$
J(v)\geq \frac{1}{2} \nor{L^* v}X^2+ (\eps -  \nor{f}Y )\nor{v}Y.
$$
Consequently, if $\nor{f}Y \leq \eps$   the functional $J$  is nonnegative, obtaining its minimum for $\tilde v= 0$. 
}

\plav{
Otherwise, the functional $J$  attains negative values as well. Indeed, take a sequence $v_n=f/n$ and calculate  
\begin{equation*}
J(v_n)= \frac{1}{2n^2} \nor{L^* f}X^2+ \frac{1}{n}  (\eps\nor{f}Y -  \nor{f}Y^2 ).
\end{equation*}
For $n$ large enough, the right hand side of the last expression is dominated by the  term $\nor{f}Y (\eps- \nor{f}Y)/n<0$, and results in negative values of $J(v_n)$.
Consequently  the minimum $\tilde v$ differs from 0. }

\plav{
As the functional $J$ is differentiable apart from the origin, its differential at $\tilde v$ equals zero, which results in the Euler-Lagrange equation \eqref{Eul-Lagr}.
\hfill $\Box$
}

\plav{
The  case $\nor{f}Y \leq \eps$ is a trivial one, which we exclude from further analysis.
}
\plav{
In the non-trivial case ($\nor{f}Y \leq \eps$)}  the Euler-Lagrange equation implies that $\tilde u = L^* \tilde v$ is an eligible solution to the problem \eqref{appr_sol} which brings the system to the boundary of the target ball around $f$. 
\plav{Moreover, as already stated above,} it is also the solution of the minimal norm among all   $u\in X$ satisfying \eqref{appr_sol}.

 In this way, finding optimal approximative solution to \eqref{appr_sol} is equivalent to finding the minimizer of the functional $J$ given by \eqref{funct_J}, which, in turn, is equivalent to solving the corresponding   Euler-Lagrange equation  \eqref{Eul-Lagr}. For these reasons, from now on, we shall concentrate on finding efficient methods for solving the latter equation. 

\begin{remark}
It is interesting to note that  problem \eqref{appr_sol} also allows for the finite dimensional solvability. More precisely, one can show that  for any finite-dimensional subspace \plav{$E\subseteq Y$} and any target $f\in E$ there exists an approximative solution $u_E$ to  \eqref{appr_sol} such that 
\begin{equation}
\label{fd_sol}
P_E (L\plav{u_E}) = P_E(f),
\end{equation}
where $P_E$ denotes the orthogonal projection to  \plav{the space $E$. similarly, by $P_{E^\bot}=I_Y-P_E$ we denote the projection to the orhogonal compliment of $E$. 
}

Such result requires analysis of the functional 
\begin{equation*}
  \label{funct_JE}
J_E(v)=\frac{1}{2} \nor{L^* v}X^2+ \eps \nor{\plav{P_{E^\bot}(v)}}Y - \skp{f}{v}Y, \quad v\in Y.
\end{equation*}
Its coercivity is proven in the manner similar to the proof of Lemma \ref{J-coerc}, which ensures existence of its \plav{unique}  minimizer $\tilde v_E$. Following the ideas of \cite{Z97}, one can show that $\tilde u_E= L^* \tilde v_E$ is the solution that satisfies both \eqref{appr_sol}  and  \eqref{fd_sol}.

\plav{
Indeed, take an arbitrary $z\in Y$. Then for any scalar $\lambda>0$ we have
\begin{equation*}
J_E( \tilde v_E) < J_E( \tilde v_E + \lambda z).
\end{equation*}
By expanding the terms in the last relation it follows
\begin{equation*}
\begin{aligned}
\eps \nor{P_{E^\bot}(\tilde v_E)}Y < &\frac{\lambda^2}2 \nor{L^*z}X^2 +\lambda  \skp{L^* \tilde v_E}{L^*z}X  \\
& +\eps \nor{P_{E^\bot}(\tilde v_E+\lambda z)}Y- \lambda\skp{f}{z}Y.
\end{aligned}
\end{equation*}
Consequently, by taking $\lambda>0$ we get
\begin{equation*}
\begin{aligned}
\skp{f}{z}Y- \skp{L^* \tilde v_E}{L^*z}X \leq & \eps\liminf_{\lambda\to 0 ^+}\frac{\nor{P_{E^\bot}(\tilde v_E+\lambda z)}Y-\nor{P_{E^\bot}(\tilde v_E)}Y}\lambda\\
&= \eps \nor{P_{E^\bot} (z)}Y.
\end{aligned}
\end{equation*}
Similarly, by taking $\lambda<0$, the same procedure implies
\begin{equation*}
\Big|\skp{f}{z}Y- \skp{L^* \tilde v_E}{L^*z}X\Big| \leq \eps \nor{P_{E^\bot} (z)}Y.
\end{equation*}
From here, by taking an arbitrary $z\in E$ the relation \eqref{fd_sol}  follows. 
}
\plav{The approximation constraint  \eqref{appr_sol} is then obtained by varying $z$ in the orthogonal complement of $E$.
}

\end{remark}

The structure of the optimal solution can be easily seen and analysed  in terms of Fourier coefficients. To this effect, let us suppose that the operator $\Lambda=  L L^* $ is diagonalisable, and denote by $\lambda_k$ \plav{the corresponding eigenvalues}. Then the  next formula  follows directly \plav{from \eqref{Eul-Lagr}}
$$
\tilde v_k = {f_k \over \lambda_k + \eps/\nor{\tilde v}Y}=  {f_k - \eps_k \over \lambda_k },
$$
where $ \eps_k= \eps \tilde v_k /\nor{\tilde v}Y$, while $\tilde v_k $  and $f_k$ denote the $k$-th  Fourier coefficient of vectors $\tilde v$  and  $f$, respectively.
Although the formula is not spectrally decomposed (expression for \plav{$\tilde v_k$} requires knowledge of $\nor{\tilde v}Y$, i.e. of all the Fourier coefficients), lack of particular frequency in the target $f$ implies that \plav{the corresponding} Fourier coefficient  of the solution $\tilde v$ vanishes as well.
\plav{In} particular, if $f$ belongs to some finite-dimensional subspace spanned by finite number of eigenvectors of  $\Lambda$, so does the solution $\tilde v$. 

Furthermore, the expression for the solution is almost explicit, up to a scalar $\nor{\tilde v}Y$. 
Although it might look surprisingly at first, this is in accordance with known results for the optimal control  problems for parabolic equations (e.g. \cite{LM21,GLNT21}). Note that such problems can be written in the form \eqref{sys_eq} with $L$ being the control to (the final) state operator. If the control acts through initial data, then this  can be considered as an inverse problem (of initial source identification). It is an important, but also numerically challenging issue due to the dissipative nature of such equations.

If the eigendecomposition of the operator $\Lambda$ is available, then the numerical procedure of calculating the optimal solution can be reduced to solving the equation for the unknown scalar by some suitable method. 
 For general  operators, with variable coefficients and/or acting on irregular domains such   decomposition is  not available or hard to construct. Another  \plav{numerical approach} should be used in that case, which \plav{usually employs some} iterative method (e.g., conjugate gradient).

\section{Parameter dependence }
\label{sect_main}
\subsection{The problems setting and characterisation of the solutions}

In the sequel we want to analyse  a family of problems of the type
\begin{equation}
  \label{sys_par}
 L_\nu u_\nu=f_\nu,
\end{equation}
where $\nu$ is the parameter ranging over a compact, \plav{connected}  set $\lN \subset {\R}^d$, $d\ge 1$. 

The goal is to propose an efficient method for finding an optimal approximative solutions to \eqref{sys_par} for a large number of parameters. To this effect, we make the following 
assumptions. 

\begin{itemize}
\item[(A1)]  
$ L_\nu $ belongs to $L(X,Y)$ for every parameter, where $X$ and $Y$ are $\nu$-independent Hilbert spaces. 
\item[(A2)] 
The associated  adjoint operator $L^*_\nu \in L(Y, X)$ is injective  for every $\nu \in \lN$ and the family of operators $  \Lambda_\nu= L_\nu L_\nu^*$ are uniformly bounded in $L(Y)$ from below by a positive self-adjoint operator  $\Lambda$, i.e.
\begin{equation*}  
\label{D_unif}
\Lambda  \leq \Lambda_\nu , \quad \nu \in \lN. 
 \end{equation*}
(The last inequality means that \plav{$\skp{(\Lambda_\nu - \Lambda) v}{v}Y\geq 0 $} for every $v\in Y$.)
\item[(A3)] The right hand side vectors  $f_\nu$ belong to a precompact subset of  $Y$. In addition, we assume uniform boundedness from below, i.e.  there exist a positive constant $f_- > \eps$ such that 
$\nor{f_\nu}Y \geq f_-.$ The aim of the last assumption is  to exclude trivial solutions and singularities in the Euler-Lagrange equation (cf. Remark \ref{triv_sol}). 
\item[(A4)] The mappings $\nu \to (L_\nu, f_\nu) \in L(X,Y) \times Y$ are analytic. 
\end{itemize}

According to \plav{the Hahn-Banach theorem}, the injectivity assumption on $L_\nu^*$ implies that the image of $L_\nu$ is dense in $Y$. Thus, for a general target $f_\nu$, the problem \eqref{sys_par} is not well posed. 
Therefore, as it was discussed in the previous section, we replace it with the following constrained optimization problem
\begin{equation}
  \label{appr_prob}
\min_{u\in X} \left\{ \nor uX \big| \nor{L_\nu u-f_\nu}Y \leq \eps\right\}.
\end{equation}
Supposing that $ \nor{f_\nu}Y $ is strictly larger than $\eps$ for every parameter value, we know that the optimal solution $\tilde u_\nu $ equals  $ L_\nu^* \tilde v_\nu$, where $\tilde v_\nu$ is the
 the solution to  the  corresponding   Euler-Lagrange equation
\begin{equation}
  \label{Eul-Lagr_par}
  L_\nu L_\nu^* \tilde v_\nu+ \eps {\tilde v_\nu \over \nor{\tilde v_\nu}Y} - f_\nu=0.
  \end{equation}
  The set of solutions we denote by  $\tilde v(\lN)=\{\tilde v_\nu, \nu\in \lN\}$. 
Equivalently, the solution $ \tilde v_\nu$ can be characterised as the minimiser of the parameter-dependent  functional $J_\nu$ defined by 
\eqref{funct_J}, with $L$ and $f$ being replaced by their parameter-dependent counterparts $L_\nu$ and $f_\nu$, respectively. 
The next result provides the boundedness of the set $\tilde v(\lN)$  in $Y$.

\begin{lemma}
\label{bound_sol}
The solutions $ \tilde v_\nu$  are uniformly  bounded in $Y$. More precisely,  there exist positive constants $v_+, v_-$ such that for every  $\nu\in\lN$ we have $v_- \leq \nor{\tilde v_\nu}Y \leq v_+.
$
\end{lemma}
{\bf Proof:} 
From the   Euler-Lagrange equation \eqref{Eul-Lagr_par} we get 
$$
\nor{L_\nu L_\nu^* \tilde v_\nu}Y  \geq \nor{f_\nu}Y -  \eps,
$$
implying the required  lower bound with $v_-= (f_- - \eps)/L_+$, where $f_- >\eps$ is the bound from the assumption (A3), and $L_+=\max_\nu \nor{L_\nu L_\nu^*}{L(Y)}$. The latter number exists as $\nu\to L_\nu$ is a smooth mapping on a compact set. 

In order to obtain  the upper \plav{bound, let} us assume the contrary. We suppose there \plav{exists} a sequence $(\nu_n)$ such that $\nor{\tilde v_{\nu_n}}Y \to \infty$. 
By multiplying  the   corresponding Euler-Lagrange equations with 
$v_{\nu_n}^0= \tilde v_{\nu_n} / \nor{\tilde v_{\nu_n}}Y$ one gets
\begin{equation}
\label{el_seq}
  \nor{v_{\nu_n}}Y \nor{L_{\nu_n}^* v_{\nu_n}^0}X^2+ \eps  - \skp{f_{\nu_n}}{v_{\nu_n}^0}Y=0.
\end{equation}
Divergence of  the sequence $ (\nor{v_{\nu_n}}Y)$ implies that $\nor{L_{\nu_n}^* v_{\nu_n}^0}X   \to 0 $. 
By exploring the assumption (A2) and the sandwich theorem we obtain that \plav{$\skp{\Lambda v_{\nu_n}^0 }{ v_{\nu_n}^0}Y \to   0 $.}

On the other hand,  denoting by $v^0$ a weak limit of $v_{\nu_n}^0$ (up to  a subsequence), we have that 
$$
\plav{\skp{\Lambda v^0}{v^0}Y }= \nor{L^* v^0}X^2 \leq \liminf_n \nor{L^* v_{\nu_n}^0}X^2  =  \lim_n \plav{\skp{\Lambda v_{\nu_n}^0}{v_{\nu_n}^0}Y }= 0.
$$
The positivity of $\Lambda$ implies $v^0=0$.

Going back to \eqref{el_seq} and exploring \plav{the relative precompactness of the set $\{f_\nu | \nu\in \lN\}$,} we obtain
$$
0> \eps  - \skp{f_{\nu_n}}{v_{\nu_n}^0}Y \to \eps ,
$$
which provides the required contradiction. 
\hfill $\Box$

In the next step we explore the smoothness properties of the solution mapping $\nu \to  \tilde v_\nu$. The smooth dependence of the solutions on the parameter is a bit delicate due to the non-smooth term in the Euler-Lagrange equation. However, the solution vanishes only when the norm of the target is small, i.e. when   $\nor{f_\nu}Y <\eps$, which is the case we excluded from the analysis through the assumption (A3). Thus, in practice there is no singularity in the
parameter-dependent Euler-Lagrange equation \eqref{Eul-Lagr_par}. 

This observation plays a key role when analyzing the smooth dependence of the solutions on the parameter. Indeed, we claim
that the solutions of \eqref{Eul-Lagr_par} depend as smoothly on the parameter as the operators $L_{\nu}$ and vectors $f_\nu$ do. To see this, we apply the implicit function theorem
(IFT). In order to do it successfully we have to analyse the nonlinear function
$$
F(v, \nu):=  L_\nu L_\nu^* v+ \eps  {v  \over \nor vY}.
 $$
Writing the solution $\tilde v_\nu$ implicitly in terms of $\nu$ needs \plav{the  computation of} the derivative of $F$ with respect to $v$ and checking that it is invertible. This
is indeed the case since
$$
\partial_v F(v, \nu):=  L_\nu L_\nu^*+ \eps \left( {I   \over \nor vY} -   {v\otimes v   \over \nor vY^2} \right)=  L_\nu L_\nu^*+ {\eps     \over \nor vY} (I_Y - P_v),
 $$
where $I_Y$ is the identity, while \plav{$P_v$ is the orthogonal projection operator on the space spanned by $v$}. 

We emphasize, despite the apparently singular term arising in the denominator, that singularity does not actually occur since we are always working with non-trivial solutions $v\not= 0$. 
The IFT can then be applied without difficulty since the resulting operator $L_\nu L_\nu^*+ {\eps     \over \nor \cdot Y}  (I_Y - P_v)$  is invertible. This can be easily seen by the Lax-Milgram lemma.

The direct application of the Lax-Milgram Lemma requires  the coercivity of $\partial_v F(v, \nu)= L_\nu L_\nu^*+ {\eps     \over \nor vY}  (I_Y - P_v)$. However, in our setting in which only an approximate solution to the original system exists, $L_\nu$ can not be guaranteed to be \plav{coercive and neither} the projection $ I_Y - P_v$ is coercive. In fact, it degenerates along the direction $v$ by the very definition of the projection. Thus it is unclear whether the addition of the term $I_Y -P_v$ to the operator $\Lambda_\nu$ suffices to guarantee the invertibility of $\partial_v F(v, \nu)$.  The positive answer is provided by the following lemma.

\begin{lemma}
\label{coerc_lemma}
Assume $Q=Q_1+Q_2$ is a bounded linear operator on a Hilbert space $H$, with $Q_1$ being a positive linear operator, while $Q_2$ is a self adjoint \plav{operator}   which vanishes on a finite dimensional subspace $V_1$ and is coercive on  $V_2=V_1^\perp$. Then the operator $Q$ is coercive on $H$.
\end{lemma}
{\bf Proof:} 
In order to prove the lemma we have to show that there exists a constant $c>0$ such that for every $v\in H$ it holds
$$
\plav{\skp{Qv}vH} \geq c \nor vH^2.
$$
We argue by contradiction. Assuming the contrary there exists a sequence $(v_n)$ such that 
$$
\plav{\skp{Qv_n}{v_n}H} \leq {1 \over n}\nor{v_n}H^2.
$$
Let us introduce the decomposition $v_n= v_{n,1}+v_{n,2}$ ,where $v_{n,i}\in V_i, i=1,2$. By using the assumptions of the lemma, it follows
$$
\plav{\skp{Qv_n}{v_n}H} = \plav{\skp{Q_1 v_n}{v_n}H}  + \plav{\skp{Q_2 v_{n,2}}{ v_{n,2}}H} \leq {1 \over n}\nor{v_n}H^2.
$$
Dividing the last \plav{equation by} $\nor{v_n}H^2$ and exploring the coercivity of $Q_2$ on the subspace $V_2$ we obtain 
\begin{equation}
\label{Q1}
\plav{\skp{Q_1 v_n^0}{v_n^0 }H}+ c_2 {\nor{v_{n,2}}H^2 \over \nor{v_n}H^2}\leq {1 \over n}\nor{v_n}H^2,
\end{equation}
where $v_n^0$ stands for the normalized  vector  $v_n / \nor{v_n}H$, while $c_2$ is the coercivity constant of $Q_2$ (on $V_2$). 

The positivity of $Q_1$ implies ${v_{n,2}\over \nor{v_n}H} \to 0$ strongly in $H$. As $V_1$ is finite dimensional, we obtain the strong convergence of the whole normalized sequence
$ v_n^0 \to v_1^0 \in V_1\setminus\{0\}$. 

Finally, by passing to the limit in \eqref{Q1} we obtain $\plav{\skp{Q_1 v_1^0}{v_1^0}H} =0 $ which contradicts the positivity of $Q_1$. 
\hfill $\Box$

As discussed above, the last lemma allows application of IFT. \plav{In particular, we employ} its analytic version (e.g. \cite[\S Appendix B]{PT87}), by which we preserve the smoothness imposed by assumption (A4). We summarize the results of this discussion in the following proposition. 
\begin{proposition}
\label{regul_sol}
The solutions to the Euler--Lagrange equation \eqref{Eul-Lagr_par} preserve the regularity imposed on the mapping $\nu \to (L_\nu, f_\nu)$. In particular, if the latter is   analytic, the same holds for the solution mapping. 
\end{proposition}

\subsection{Greedy approach and the Tychonoff regularization}
The aim of this section is to develop an efficient numerical algorithm for reconstructing an arbitrary  element of $\tilde v(\lN)$ corresponding to some given parameter value. In order to accomplish this task we rely on  greedy algorithms which were introduced and analysed through the last two decades  in the context of parametric PDEs. They serve as on of the  most popular tools for construction of  reduced basis  (cf. \cite{LL21} and the references therein). 

The objective of the greedy approach is to approximate a compact set $\lK$ (e.g.  a family of solutions to parameter dependent problems) in a Banach space $Y$ by a linear  subspace  $V_n$ of (small) dimension $n$. The selection of its basis vectors is done gradually in the offline phase of the algorithm. Once the basis is known, the approximation for an arbitrary given element in \plav{$\lK$ } is computed in the online phase. Usually, the computational effort for the offline routine is much higher than for the online one, but it is performed only once.

A greedy approximation is optimal, where the optimality is to be understood  in the sense of the Kolmogorov widths. The Kolmogorov $n-$width defines (theoretically) the best possible approximation error one could obtain by a subspace
in $Y$ of a fixed dimension $n$. The results of \cite{Bin11,DPW13} show that  the greedy approximation errors decay asymptotically with 
the same rate (exponential or polynomial) as the Kolmogorov widths. 

Furthermore, the Kolmogorov widths are preserved under smooth (analytic) mappings (\cite{CD}). In such a way one can a-priori estimate approximation performance of a subspace constructed by a greedy method, by estimating Kolmogorov widths of a set of admissible parameters that generate a set of interest. 
However, it is important to emphasize that the very implementation of a greedy procedure does not require analytic smoothness.  Only, in that case, we lack the a-priori estimates on the  approximation
errors.
\plav{In particular, this might happen if the parameter set $\mathcal N$ consists of a finite or a countable number of elements, which is not the case studied in this article.  }

\plav{
In the development of greedy algorithms one of the main steps is the construction of some  surrogate function which enables us to calculate the distance between unknown terms, or terms that are } in general hard to calculate. Usually such a surrogate is provided by using some kind of residual. More, precisely, let us assume we have calculated $\tilde v_{\nu_1}$ for some parameter value $\nu_1$. We would like to check if we can use it to approximate some other optimal vector $\tilde v_{\nu}$. To this effect we plug  $\tilde v_{\nu_1}$ \plav{into the equation} satisfied by the latter term and define the residual
$$
R_\nu  \tilde v_{\nu_1} :=   L_\nu L_\nu^* \tilde v_{\nu_1}+ \eps N  \tilde v_{\nu_1} - f_\nu,
$$
where by $N$ we denote  the normalization operator $N v  =v/ \nor{v}Y$.
If the residual turns to be zero, due to the uniqueness of the solution it follows that $\tilde v_{\nu} = \tilde v_{\nu_1}$. In general we would like the residual to measure the distance between two optimal vectors. More precisely, we require estimates of the form 
\begin{equation}
  \label{res_est}
c_- \nor{\tilde v_\nu-\tilde v_{\nu_1}}Y \leq \nor{R_\nu  \tilde v_{\nu_1}}Y  \leq   c_+ \nor{\tilde v_\nu-\tilde v_{\nu_1}}Y,
  \end{equation}
where $c_-$ and $c_+$ are $\nu$-independent positive constants.  

In order to obtain such kind of estimates let us rewrite the residual operator as
\begin{equation}
  \label{res}
R_\nu  \tilde v_{\nu_1} :=   L_\nu L_\nu^* (\tilde v_{\nu_1} - \tilde v_{\nu}) + \eps (N  \tilde v_{\nu_1} - N  \tilde v_{\nu}),
 \end{equation}
where we explored the Euler-Lagrange equation \eqref{Eul-Lagr_par}. 
The upper bound in \eqref{res_est} now follows easily. For the first summand in \eqref{res} it is a direct consequence of the boundedness assumption, while for the 
 last term in \eqref{res} it is obtained by using geometrical interpretation and the triangular inequality (cf. \cite[\S 1.2]{AL13}).

However, due to the structure of the normalization operator (which is constant along each half-ray emerging from the origin), \plav{one can easily check that the lower bound in \eqref{res_est} is equivalent to the coercivity} of the operator $ L_\nu L_\nu^*$. 
Of course, our assumptions on $L_\nu$ do not provide the required coercivity and we have to propose an alternative method. Note that here we can not use the approach applied in Lemma \ref{coerc_lemma} as the normalization operator $N$ does not vanish on any non-trivial subspace. 

In order to overcome the lack of coercivity, let us introduce a two-parame\-ter family of linear problems
\begin{equation}
  \label{2par-eq}
  L_\nu L_\nu^* \tilde w(\nu, \eta)+ \eta   \tilde w(\nu, \eta) = f_\nu.
    \end{equation}
    Here we have substituted the nonlinear normalization operator $N$ appearing in the Euler-Lagrange equation \eqref{Eul-Lagr_par} by a linear term multiplied by the  new parameter $\eta$. Let us note that the last equation coincides with the Euler-Lagrange equation obtained by  minimization of the functional 
    \begin{equation}
  \label{funct_J_reg}
J(v)=\frac{1}{2} \nor{L_\nu^* v}X^2+ \eta \nor{v}Y^2 - \skp{f_\nu}{v}Y, \quad v\in Y.
\end{equation}
The non-smooth term appearing in the original functional \eqref{funct_J} is here replaced by a quadratic one. This  improves the coercivity properties of the problem, which is a standard benefit of the Tychonoff regularization.  

\plav{
The functional appearing in \eqref{funct_J_reg} corresponds to the penalization approach for approximation problems, where one forces the solution to approach the given target by letting the penalization constant blow up. It provides a smooth functional that is easier to handle, unlike the one in \eqref{funct_J}. However, such an approach does not allow an a-priori estimate of the deviation from the target, in particular, it does not ensure that it is smaller than the given approximation error $\eps$. For this reason, the value of the parameter $\eta$ in (3.9) is not fixed, but is allowed to vary within a specific interval.
}

    \plav{ More precisely, we suppose the   introduced  parameter} ranges within the interval $\zi{\eta_-}{\eta_+} =\zi{\eps/v_+}{\eps/v_-}$, where $v_\pm$ are bounds from Lemma \ref{bound_sol}. In such a way, for each parameter $\nu$ there exists an $\eta$ from the given range such that $\eta= \eps/\nor{\tilde v_\nu}Y$, where $\tilde v_\nu$ is the corresponding solution of the Euler-Lagrange equation \eqref{Eul-Lagr_par}. 
    
    This implies that solving a family of two-parameter problems \eqref{2par-eq} for $(\nu, \eta)= \lN\times \zi{\eta_-}{\eta_+}$ will also provide the solution of the original $\eps$ problem. 
    
   In order to efficiently  treat the auxiliary problem \eqref{2par-eq} we explore greedy algorithms introduced above.     The problem is now linear, and it involves selfadjoint operators of the form $L_\nu L_\nu^* + \eta I$ which are uniformly bounded from below by $\eta_- I$. This allows one to consider the residual of the form
   $$
   R(\nu, \eta) w:=  L_\nu L_\nu^*  w+ \eta    w- f_\nu
   $$ 
   with $w$ being an arbitrary test function. Based on the above obtained bounds it directly follows 
    \begin{equation*}
  \label{res_est2}
c_-  \nor{\tilde w(\nu, \eta) -w}Y \leq  \nor{ R(\nu, \eta) w}Y \leq   c_+ \nor{\tilde w(\nu, \eta) -w}Y.
  \end{equation*}
  Furthermore, employing the same kind of arguments based on the IFT and presented in the previous subsection, the solutions $\tilde w(\nu, \eta)$ to  \eqref{2par-eq} preserve the smoothness of the mapping $\nu \to (L_\nu, f_\nu)$ at all levels. However, note that in this case the corresponding analysis is much simpler as there is no singularity in the two parameter equation  \eqref{2par-eq} and we deal with uniformly coercive operators.

  In particular if the operators $L_\nu$ and non-homogenous terms $f_\nu$ depend analytically on $\nu$, then the mapping $(\nu, \eta) \to \tilde w(\nu, \eta)$ is analytic as well.
  In this way, as discussed at the beginning of this section,  one preserves the Kolmogorov widths of the two-parameter set $ \lN\times \zi{\eta_-}{\eta_+}$ which are transferred to the manifold of solutions $\tilde W=\{\tilde w(\nu, \eta) \big| (\nu, \eta)= \lN\times \zi{\eta_-}{\eta_+} \}$ (cf. \cite{CD}).

Suppose we have \plav{performed a greedy algorithm} for the two-parameter problems \eqref{2par-eq}. It selects a finite set of parameter pairs $(\nu_i, \eta_i), i=1..N$ and returns the corresponding solutions $w_i:=\tilde w(\nu_i, \eta_i)$ which constitute a reduced basis for the manifold $\tilde W$. In other words,  for every value of $(\nu, \eta)$ there exist a set of linear coefficients $\alpha_i$ such that 
$\nor{\tilde w(\nu, \eta) - \sum \alpha_i w_i }Y< \delta$ where $\delta$ is a positive constant determined by the stopping criteria of the greedy algorithm.

    In the next step we want to go back to the original $\eps$ problem \eqref{Eul-Lagr_par}. More precisely, given an arbitrary value of the parameter $\nu$ we want to determine a set of coefficients $\alpha_i$ such that $v_\nu^\star=\sum \alpha_i w_i$ brings the system within $\eps$ distance from the target $f_\nu$. The problem is feasible, as for $\eta=\eps/\nor{\tilde v_\nu}Y$ the greedy algorithm provides  a good approximation of the solution to \eqref{2par-eq}, and consequently to the original Euler-Lagrange equation \eqref{Eul-Lagr_par}. 
    
    Therefore, let us propose the required approximation by projecting the target  $f_\nu$ on the space spanned by $L_\nu L_\nu^* w_i$. More precisely, we determine the approximation coefficients $\alpha_i$ as solutions to the system
    $$
    \sum_{i=1}^{\plav{m}} L_\nu L_\nu^* \alpha_i w_i = \plav{P_\nu^m} f_\nu,
    $$
    with \plav{$P_\nu^m$} denoting the orthogonal projection on the space spanned by vectors $L_\nu L_\nu^* w_i$, while \plav{$m$} stands for the number of parameters selected during the offline phase of the greedy procedure. The approximation of the solution $\tilde u_\nu$ to the constrained optimal control problem \eqref{appr_prob} is then given by
    $$
    \tilde u_\nu \approx   \sum_{i=1}^{\plav{m}}  L_\nu^* \alpha_i w_i.
    $$
    Such procedure results in an approximate solution that steers the system to the target $f_\nu$ as close as possible by means of the constructed reduced basis space $\{w_i, i=1..\plav{m}\}$.

\section{Unbounded operators}

The theory developed in the previous section requires the operators of interest to be bounded. However, it can be generalised and  applied to unbounded operators as well.

In this section we consider  a family of  linear unbounded operators $A_\nu$ on a Hilbert space \plav{$H$}, where, as before,  $\nu$ is the parameter ranging over a compact, \plav{connected} set $\lN \subset {\R}^d$, $d\ge 1$. 
We put the following hypothesis on the considered family. 
\begin{itemize}
\item[(H1)]  $A_\nu$ are positive, self adjoint operators uniformly bounded from below, i.e. there exists $\alpha\in \R^+$ such that $A_\nu \geq\alpha>0$ for every $\nu\in \lN$;
\item[(H2)]  The operators $A_\nu$ have common domains, i.e. there exists a subspace \plav{$D_A\subseteq H$}  such that $D(A_{\nu}) =D_A$ for every $\nu\in \lN$;
\item[(H3)] the graph norms of $A_\nu$  are  uniformly equivalent. 
\end{itemize}
\begin{remark}
\label{rem-A1}
\plav{Instead of  (H1) one can  require the operators  $A_\nu$ to be  uniformly lower (or upper) bounded by an arbitrary constant.} In order to simplify the presentation we restrict to the  case of positive definite operators.  
\end{remark}
 In addition we suppose that  $H$ is densely and compactly embedded into some Hilbert space $Y$. We pose the problem of finding (an approximative) solution to the equation
 \begin{equation*}
A_\nu x = y, \quad y\in Y.
\end{equation*}
Based on the assumption (H1) the image of  operators $A_\nu$  equals $H$, but it is only dense in $Y$. This brings us to the problem of finding the optimal approximative solution discussed in previous sections. However, the theory we developed assumes bounded operators, which $A_\nu$ are not. 
In order to overcome this gap,   we want to associate to each $A_\nu$ an operator $L_\nu: X\to Y$, where $X$ is some  parameter independent \plav{Hilbert space} which is still  to be defined.

To this effect, let us denote by $X_\nu$ the space $D_A$ equipped with the norm
\begin{equation}
\label{X1_nor}
\nor x{X_\nu}= \nor{(\beta I_H - A_\nu)x}H, 
\end{equation}
where $\beta<0$ is a scalar from the resolvent set of $A_\nu$, while $I_H$ stands for the identity on $H$. Note that the introduced space corresponds to the one denoted by $X_1$ in \cite[\S 2.10]{TW09}, and 
their norms  are uniformly equivalent to \plav{the graph norms} of $A_\nu$  (\cite[Proposition 2.10.1]{TW09}).
\plav{Based on the assumption (H3) it follows that the norms \eqref{X1_nor} are} uniformly equivalent. For this reason in the sequel we use a common notation $X$ for all spaces $X_\nu$.

Let us denote by $\tilde L_\nu$ a family of operators from $X$ to $H$ defined by
\begin{equation*}
\label{tL_def}
\tilde L_\nu x = A_\nu x.
\end{equation*}
By the definition of \plav{the space $X$}, it  is not difficult to check that the introduced operators are uniformly bounded and coercive. 
Indeed, we have
\begin{equation}
\label{L_ub}
\nor{\tilde L_\nu  }{L(X,H)}=\sup_{x\in D_A} {\nor{\tilde L_\nu x }H  \over \nor xX}
=\sup_{x\in D_A} {\nor{A_\nu x }H \over \nor{(\beta I_H - A_\nu)x}H} \leq 1,
\end{equation}
where the last inequality follows from the positivity of the operator $A_\nu$. Similarly, by exploring 
\begin{equation}
\label{L_lb}
 {\nor{\tilde L_\nu x }H  \over \nor xX}
\geq  {\nor{A_\nu x }H \over \nor{\beta x}H + \nor{A_\nu x }H }  \geq  {\alpha \over \beta+\alpha},
\end{equation}
where $\alpha$ is the lower bound from the assumption \plav{(H1)}, one obtains the lower bound on $\tilde L_\nu$. 


Thus we obtain that $\tilde L_\nu$ form a family of uniformly bounded operators in $L(X,H)$. 
Based on the assumption (H1) both the operators $A_\nu$  and  $\tilde L_\nu$ are surjective into $H$. For the corresponding adjoints we have the following characterization. 

\begin{lemma}
\label{L*}
The adjoint operators $\tilde L_\nu^*\in L(H, X)$ are of the form $\tilde L_\nu^* = A_\nu (\beta I_H - A_\nu)^{-2}  \in L(H, X)$ , satisfying the same lower and upper bounded estimates as the operators $\tilde L_\nu$ do. 
\end{lemma}
{\bf Proof:} 
By the definition of the adjoint, for $u, x \in  X$  we have
\begin{equation*}
\begin{aligned}
\skp x{\tilde L_\nu^* u}X &= \skp {\tilde L_\nu x}uH = \skp {A_\nu x}uH =  \skp x{A_\nu u}H\\
&= \skp { (\beta I_H - A_\nu)^{-1} x }{(\beta I_H - A_\nu)^{-1}A_\nu u}X.
\end{aligned}
\end{equation*}
From here we get 
$$
\tilde L_\nu^* u = A_\nu (\beta I_H - A_\nu)^{-2} u, \quad u\in X.
$$
As $X$ is dense in $H$, and the operators $A_\nu (\beta I_H - A_\nu)^{-1}$ are bounded on $H$, the last relation holds for an arbitrary $u\in H$. 
This provides the first part of the statement.

By using the obtained explicit expression for the adjoint operators, we have
$$
{\nor{\tilde L_\nu^* u }X  \over \nor uH} = {\nor{A_\nu  (\beta I_H - A_\nu)^{-1} u  }H  \over \nor uH} = {\nor{A x }H \over \nor{(\beta I_H - A)x}H} ,
$$
where $x= (\beta I_H - A_\nu)^{-1} v\in X$. The required bounds now follow from \eqref{L_ub} and  \eqref{L_lb}.  
\hfill $\Box$

Finally, in order to  put us in the context of the preceding section, we introduce  $L_\nu\in L(X, Y)$ defined as 
\begin{equation}
\label{L_dec}
L_\nu:= I  \tilde L_\nu,
\end{equation}
where $I$ stands for the inclusion operator $I: H \to Y$. As $H$ is  compactly  embedded into  $Y$, the operator $L_\nu$ is bounded, injective, but not coercive operator with the  dense image  in $Y$. Moreover, due to the uniform  coercivity of operators $\tilde L_\nu^*$
we have that 
\begin{equation*}
\label{}
L_\nu L^*_\nu \geq c I I^*.
\end{equation*}
This implies the assumption (A2) is fulfilled and consequently we fit the setting of the previous section.

\section{A specific example}
We consider a family of Dirichlet Laplacians  on $L^2(\Omega)$, where $\Omega$ is assumed to be  an open, bounded set with a smooth boundary.
More precisely, we define a family of unbounded operators on $L^2(\Omega)$ 
\begin{equation}
\label{Lapl}
-\Delta_\nu= - \dv (\mA_\nu \nabla)
\end{equation}
accompanied by Dirichlet boundary conditions. For the coefficients $\mA_\nu$ we assume that they depend smoothly on $\nu$ and they satisfy uniform boundedness and coercivity properties. More precisely, we assume $\mA_\nu\in C^1(\Omega)$ are self adjoint matrix functions  such that 
$$
\mA_- (x) \leq \mA_\nu(x)  \leq \mA_+(x) , \quad \nu\in \lN, x \in \Omega,
$$
for some bounded and coercive matrix functions $\mA_\pm\geq \alpha>0$. By using the Poincare inequality this immediately applies the uniform coercivity of the considered Laplacians, i.e. $-\Delta_\nu \geq \alpha$. 

The domain of the Laplacian $-\Delta_\nu$ is parameter independent and \plav{coincides} with $X=H^2(\Omega) \cap H^1_0(\Omega)$ for every $\nu$. 
Furthermore, by standard elliptic regularity results (e.g. \cite[\S 6.3]{E10}), their graph norms are (uniformly) equivalent to $H^2$ norm. Consequently, the hypo\-thesis (H1)-(H3) from the previous section are satisfied. 

For the target space we take $Y=  H^{-1}(\Omega)$ for which we have the dense and compact embedding $L^2(\Omega) \hookrightarrow H^{-1}(\Omega)$ (of course, any other space $H^s(\Omega)$, with $s<0$ , will be appropriate in this context). Consequently, we introduce a sequence of bounded operators $L_\nu \in L(X, Y)$ defined by
$$
L_\nu u =  - \dv (\mA_\nu \nabla u).
$$
The properties of the introduced operators can be examined \plav{through their} matrix representation. To this effect, we explore the spectral decomposition of the Laplacian operator \eqref{Lapl}.
In particular, there  exists an orthonormal basis in  $L^2(\Omega)$ consisting of eigenfunctions $\psi_{i, \nu}$ of $-\Delta_\nu$ such that
$$
-\Delta_\nu \sim \diag(\lambda_{1, \nu}, \lambda_{2, \nu}, \ldots),
$$
where $(\lambda_{i, \nu})$ is a sequence of (positive) eigenvalues diverging to infinity. 

Then it is not difficult to check that the matrix representation of the associated operator $L_\nu$ in the pair of basis 
$\hat \psi_{i, \nu}: = \psi_{i, \nu} /\lambda_{i, \nu}\in H^2(\Omega) \cap H^1_0(\Omega)$ and $\tilde \psi_{i, \nu}: = \sqrt{\lambda_{i, \nu}}\psi_{i, \nu} \in H^{-1}(\Omega)$ has the form
$$
L_\nu \sim \diag({1 \over\sqrt{\lambda_{1, \nu}}}, {1 \over\sqrt{\lambda_{2, \nu}}}, \ldots).
$$
Similarly, the same matrix representation form, in the reverse pairs of basis $\{\tilde \psi_{i, \nu}\}$ and $\{\hat \psi_{i, \nu}\}$, also holds for the adjoint operator $L_\nu^*$. Due to the properties of eigenvalue sequence, the operator $L_\nu^*$ is injective, but not coercive. This is equivalent to the statement that the image of the Laplacian is (only) dense in  $H^{-1}(\Omega)$.

Having associated to Laplacians \eqref{Lapl} a sequence of bounded operators $L_\nu$ with dense images, we can explore the approach developed  in Section \ref{sect_main} for finding the optimal approximative solution.
To this effect, we assume $f_\nu: \lN \to H^{-1}(\Omega)$ is a smooth function with a precompact image 
and we consider a sequence of problems
\begin{equation}
\label{prob_lapl}
 - \dv (\mA_\nu \nabla) u = f_\nu,
\end{equation}
with solutions $u$ searched within  the domain of the Laplacian, i.e. in  $X=H^2(\Omega) \cap H^1_0(\Omega)$. As $L_\nu(X)$ is only dense in $H^{-1}(\Omega)$, the problem in general does not admit a solution, and we  relax  it by considering an approximative one 
of the form \eqref{appr_prob}. 
The optimal  solution is thus of the form   $\tilde u_\nu=  L_\nu^* \tilde v_\nu$, where  $ \tilde v_\nu$ is the solution of the  corresponding   Euler-Lagrange equation
\eqref{Eul-Lagr_par}.  In the next lemma we provide the explicit expression of the adjoint  operator $L_\nu^*$.

\begin{lemma}
\label{L*_ex}
  The adjoint operator equals to $L_\nu^*=(-\Delta)^{-2} (-\Delta_\nu)(-\Delta)^{-1} $, where $-\Delta$ is the Laplacian with constant coefficients $\mA= \mI$. 
\end{lemma}
{\bf Proof:} 
\plav{Let us write} the operator $L_\nu$ in the form \eqref{L_dec}, i.e. $L_\nu = I  \tilde L_\nu$. Here $I$ stands for the inclusion operator $I: L^2(\Omega) \hookrightarrow H^{-1}(\Omega)$, while $\tilde L_\nu$ is the operator from $H^2(\Omega) \cap H^1_0(\Omega)$ to $L^2(\Omega) $ defined by $\tilde L_\nu u =  - \dv (\mA_\nu \nabla u)$. 

Then for calculating the adjoint we use the relation $L_\nu^* = \tilde L^*_\nu I^*$. Similarly as in Lemma \ref{L*} one obtains that $\tilde L^*_\nu= (-\Delta)^{-2} (-\Delta_\nu)$. Thus it remains to express the adjoint of the inclusion. To this effect let us note
\begin{equation*}
\begin{aligned}
\skp u{I^* v}{L^2}&= \skp {I u}v{H^{-1}}  = \skp{(-\Delta)^{-1/2}u }{(-\Delta)^{-1/2} v}{L^2} ,
\end{aligned}
\end{equation*}
from where we get $I^* = (-\Delta)^{-1}$, which completes the proof.
\hfill $\Box$

Based on the last lemma, the Euler-Lagrange equation \eqref{Eul-Lagr_par} for the problem \eqref{prob_lapl} can be rewritten as
\begin{equation}
\label{EL_lapl}
(-\Delta_\nu) (-\Delta)^{-2} (-\Delta_\nu)(-\Delta)^{-1} \tilde v_\nu+ \eps N  \tilde v_\nu - f_\nu=0,
  \end{equation}
  where, as before,  $N$ stands for the normalization operator $N v  =v/ \nor{v}{H^{-1}}$. 

   Of course, in numerical calculations one should employ some standard discretisation procedure (based on finite differences or finite elements), and reduce the last equation to a finite dimensional, algebraic one.
\plav{   The numerical resolution using gradient methods for a specific realization of the model has been previously developed in \cite{CGL94}. Once the duality in Propositions 2.1-2.3 has been adopted in the computational method, its implementation does not differ significantly from the classical problems, that in our case would correspond to operators $L$ which, instead of being dense, are assumed to have a full range.
}

   However, as we have explained, our goal is not to solve the equation \eqref{EL_lapl} for every value of the parameter. Instead, we employ a greedy procedure which allows us to obtain a reduced basis by which we can approximate any $ \tilde v_\nu$ with some a-priori given precision. The reduced basis  is constructed during the offline phase, and consists of solutions to two-parameter problem \eqref{2par-eq} for some carefully 
   selected parameters' values. 
   
 \plav{  To the best of authors' knowledge, so far the greedy methods were addressed only to the exact solutions. On the other hand, the approximate solution problem is relevant and natural in several contexts such as image processing, control of time-dependent irreversible processes, etc. Our contribution is the analog of the greedy theory for exact solutions in the context of the approximate one (e.g. \cite{CD15} ).
 }
 
 \plav{ 
   The complete analysis of the computational cost is beyond the scope of this paper, but the well-known conclusions for the exact solution problem apply in this case as well (cf. \cite[\S 5]{LZ16}). In particular, the implementation requires extensive offline work, which is the most expensive part of the algorithm. Once this is done, the method is optimal since it leads to sharp approximation rates.
 }
 
 \plav{ 
The cost of the online part of the algorithm is of order $m\, C$, where $C$ is the cost of applying (a finite dimensional approximation of ) the operator $L_\nu$ to an arbitrary vector, while $m$ is the number of parameters selected during the offline phase. Consequently, the cost reduction obtained by choosing the greedy algorithm depends linearly on the ratio between the number of selected snapshots  and the system dimension.
 }
 
 \plav{ 
Practical application of the method is justified when the approximate equation needs to be solved for a wide class of parameter-dependent problems.
}

\section{Conclusion}

In this article we develop a procedure for handling a class of parameter dependent,  ill-posed problems, that, in general, do not \plav{allow    exact solution}. The original problem is relaxed by considering corresponding approximate problems, whose optimal solutions are well defined, where the optimality is determined through the minimal norm requirement. 

The procedure is based upon reduced basis methods, in particular upon greedy algorithms, by which one constructs a reduced basis during  the offline phase. An approximation of the solution for a specific parameter is constructed in the online phase as a suitable linear combination of reduced basis vectors. In order to provide a-priori estimates for the algorithm, a Tychonoff-kind regularization is applied, which adds an additional parameter to the model. 


The theory is developed in a rather general theoretical framework, which allows its application to different  kinds of  problems.  
As a specific example we consider a  family of ill-posed elliptic problems. The required general assumptions in this case translate to rather natural uniform lower and upper  bounds on coefficients of the considered operators. 

Other potential applications  would include approximate controllability, inverse problems related to high dissipative systems, like heat equations, deconvolution in image processing etc. Beside identifying the class of suitable operators for each of these problems, it would also be  interesting to perform corresponding numerical simulations and verify the  efficiency of the method on particular examples. 


\section*{Acknowledgements}

This research was done while the first author visited Chair of Dynamics, Control and Numerics (Alexander von Humboldt Professorship) at  Friedrich-Alexander-Universität Erlangen-Nürnberg, with the support of the DAAD (Research Stays for University Academics and Scientists, 2021 programme) and Alexander von Humboldt-Professorship. 

This project has received funding from the European Research Council (ERC) under the European Union’s Horizon 2020 research and innovation programme (grant agreement NO: 694126-DyCon), the Alexander von Humboldt-Professorship program, the European Unions Horizon 2020 research and innovation programme under the Marie Sklodowska-Curie grant agreement No.765579-ConFlex, the Transregio 154 Project ‘‘Mathematical Modelling, Simulation and Optimization Using the Example of Gas Networks’’, project C08, of the German DFG, 
 the Grant MTM2017-92996-C2-1-R COSNET of MINECO (Spain), by the Elkartek grant KK-2020/00091 CONVADP of the Basque government and by the Air Force Office of Scientific Research (AFOSR) under Award NO: FA9550-18-1-0242. 

The authors acknowledge E. Trélat  for \plav{his interesting comments} that have improved the final version of the manuscript.

\end{document}